
\documentstyle{amsppt}
\baselineskip18pt
\magnification=\magstep1
\pagewidth{30pc}
\pageheight{45pc}
\hyphenation{co-deter-min-ant co-deter-min-ants pa-ra-met-rised
pre-print pro-pa-gat-ing pro-pa-gate
fel-low-ship Cox-et-er dis-trib-ut-ive}
\def\leaderfill{\leaders\hbox to 1em{\hss.\hss}\hfill}

\

\def\a{{\alpha}}
\def\be{{\beta}}
\def\g{{\gamma}}

\def\bb{{\bold b}}
\def\bc{{\bold c}}
\def\bd{{\bold d}}

\def\bi{{\bold i}}
\def\bj{{\bold j}}

\def\bp{{\bold p}}
\def\bq{{\bold q}}

\def\bu{{\bold u}}
\def\bv{{\bold v}}

\def\brr{{\bar r}}
\def\b0{\text{\bf 0}}

\def\boxit#1{\vbox{\hrule\hbox{\vrule \kern3pt
\vbox{\kern3pt\hbox{#1}\kern3pt}\kern3pt\vrule}\hrule}}
\def\rabbit{\vbox{\hbox{\kern0pt
\vbox{\kern0pt{\hbox{---}}\kern3.5pt}}}}

\def\tableau#1{
        \hbox {
                \hskip -10pt plus0pt minus0pt
                \raise\baselineskip\hbox{
                \offinterlineskip
                \hbox{#1}}
                \hskip0.25em
        }
}

\def\tabCol#1{
\hbox{\vtop{\hrule
\halign{\strut\vrule\hskip0.5em##\hskip0.5em\hfill\vrule\cr\lower0pt
\hbox\bgroup$#1$\egroup \cr}
\hrule
} } \hskip -10.5pt plus0pt minus0pt}

\def\CR{
        $\egroup\cr
        \noalign{\hrule}
        \lower0pt\hbox\bgroup$
}



\def\blank#1#2{
\hbox to #1{\hfill \vbox to #2{\vfill}}
}


\def\strut{\vrule height10pt depth5pt width0pt}

\topmatter
\title Freely braided elements in Coxeter groups, II
\endtitle

\author R.M. Green and J. Losonczy \endauthor
\affil Department of Mathematics \\ University of Colorado \\
Campus Box 395 \\ Boulder, CO  80309-0395 \\ USA \\ {\it  E-mail:}
rmg\@euclid.colorado.edu \\
\newline
Department of Mathematics\\ Long Island University\\ Brookville,
NY  11548\\ USA\\ {\it  E-mail:} losonczy\@e-math.ams.org\\
\newline
\endaffil

\abstract  We continue the study of freely braided elements of
simply laced Coxeter groups, which we introduced in a previous
work. A known upper bound for the number of commutation classes of
reduced expressions for an element of a simply laced Coxeter group
is shown to be achieved only when the element is freely braided;
this establishes the converse direction of a previous result. It
is also shown that a simply laced Coxeter group has finitely many
freely braided elements if and only if it has finitely many fully
commutative elements. \endabstract

\subjclass 20F55\endsubjclass

\keywords braid relation, commutation class, Coxeter group, root
sequence \endkeywords

\endtopmatter

\centerline{\bf To appear in Advances in Applied Mathematics}

\head Introduction \endhead

In \cite{{\bf 5}} we defined, for an arbitrary simply laced
Coxeter group, a subset of ``freely braided elements''.  Such
elements include the fully commutative elements of Stembridge
\cite{{\bf 8}} as a particular case. The idea behind the
definition is that although it may be necessary to use long braid
relations in order to pass between two reduced expressions for a
freely braided element, the necessary long braid relations in a
certain sense do not interfere with one another.

Every reduced expression for a Coxeter group element $w$
determines a total ordering of the set of positive roots made
negative by $w$. These totally ordered sets are known as ``root
sequences''. If a root sequence for an element $w$ of a simply
laced Coxeter group contains a consecutive subsequence of the form
$\a, \a + \be, \be$, then we refer to the set of these roots as a
``contractible inversion triple" of $w$. A group element is said
to be ``freely braided" if its contractible inversion triples are
pairwise disjoint.

Let $N(w)$ denote the number of contractible inversion triples of
$w$. It was shown in \cite{{\bf 5}} that the number of commutation
classes (short braid equivalence classes of reduced expressions)
of $w$ is bounded above by $2^{N(w)}$, and that this bound is
achieved if $w$ is freely braided. In this paper, we prove that
the bound can be achieved only for freely braided $w$ (Theorem
2.2.1). This was previously known in the type $A$ setting
\cite{{\bf 5}, Theorem 5.2.1}, but the argument given here has the
advantage of being conceptual as well as more general.

The second main result of this paper is a classification of the
simply laced Coxeter groups having only finitely many freely
braided elements. Note that it is possible for such a group to be
infinite.  It turns out (see the discussion following Theorem
3.3.3) that this question has the same answer as a classification
question previously answered by others \cite{{\bf 2}, {\bf 4},
{\bf 8}}: we will show in Theorem 3.3.3 that a simply laced
Coxeter group has finitely many freely braided elements if and
only if it has finitely many fully commutative elements.  One
implication here is easy, but the converse requires some effort.
Our proof of Theorem 3.3.3 does not rely on a case analysis based
on any classification result.

\head 1. Preliminaries \endhead

\subhead 1.1 Basic terminology and notation \endsubhead

Let $W$ be a simply laced Coxeter group with distinguished
generators $S = \{s_i : i\in I\}$ and Coxeter matrix
$(m_{ij})_{i,j\in I}$. For the basic facts concerning Coxeter
groups, see \cite{{\bf 1}} or \cite{{\bf 6}}. Denote by $I^*$ the
free monoid on $I$. We call the elements of $I$ {\it letters} and
those of $I^*$ {\it words}. The {\it length} of a word is the
number of factors required to write the word as a product of
letters. Let $\phi : I^* \longrightarrow W$ be the surjective
morphism of monoid structures satisfying $\phi(i) = s_i$ for all
$i\in I$. A word $\bi \in I^*$ is said to {\it represent} its
image $w=\phi(\bi)\in W$; furthermore, if the length of $\bi$ is
minimal among the lengths of all the words that represent $w$,
then we call $\bi$ a {\it reduced expression} for $w$. The {\it
length} of $w$, denoted by $\ell(w)$, is then equal to the length
of $\bi$.

Let $V$ be a vector space over the field of real numbers with
basis $\{\a_i : i\in I\}$, and denote by $B$ the {\it Coxeter
form} on $V$ associated to $W$. This is the symmetric bilinear
form satisfying $B(\a_i,\a_j) = - \cos{\frac{\pi}{m_{ij}}}$ for
all $i,j\in I$. We view $V$ as the underlying space of a
reflection representation of $W$, determined by the equalities
$s_i \a_j = \a_j - 2B(\a_j,\a_i)\a_i$ for all $i,j\in I$. The
Coxeter form is preserved by $W$ relative to this representation.

Denote by $\Phi$ the {\it root system} of $W$, i.e., the set
$\{w\a_i : w\in W\text{ and }i\in I\}$.  Let $\Phi^+$ be the set
of all $\be \in \Phi$ such that $\be$ is expressible as a linear
combination of the $\a_i$ with nonnegative coefficients, and let
$\Phi^- = -\Phi^+$. We have $\Phi = \Phi^+ \cup \Phi^-$
(disjoint). The elements of $\Phi^+$ (respectively, $\Phi^-$) are
called {\it positive} (respectively, {\it negative}) roots. The
$\a_i$ are also referred to as {\it simple} roots. We define the
{\it height} of any root $\be$ to be the sum of the coefficients
used to express $\be$ as a linear combination of the simple roots.

Associated to each $w\in W$ is the {\it inversion set} $\Phi(w) =
\Phi^+ \cap w^{-1}(\Phi^-)$. It has $\ell(w)$ elements and
determines $w$ uniquely. Given any reduced expression
$i_1i_2\cdots i_n$ for $w$, we have $\Phi(w) = \{r_1,r_2,\dots,r_n
\}$, where $r_1 = \a_{i_n}$ and $r_l = s_{i_n}\cdots
s_{i_{n-l+2}}(\a_{i_{n-l+1}})$ for all $l\in \{2,\dots,n\}$.  Form
the sequence $\brr = (r_1,r_2,\dots,r_n)$. We call $\brr$ the {\it
root sequence} of $i_1i_2 \cdots i_n$, or a root sequence {\it
for} $w$. Notice that any initial segment of a root sequence is
also a root sequence for some element of $W$.

Let $w \in W$.  Any subset of $\Phi(w)$ of the form $\{ \a, \be,
\a + \be \}$ will be called an {\it inversion triple} of $w$. We
say that an inversion triple $T$ of $w$  is {\it contractible} if
there is a root sequence for $w$ in which the elements of $T$
appear consecutively (in some order). The number of contractible
inversion triples of $w$ will be denoted by $N(w)$.  If the
contractible inversion triples of $w$ are pairwise disjoint, then
$w$ is said to be {\it freely braided}.

\subhead 1.2 Braid moves \endsubhead

Given any $i,j\in I$ and any nonnegative integer $n$, we write
$(i,j)_n$ for the length $n$ word $iji\cdots\in I^*$. Let $\bi,
\bj \in I^*$ and let $i,j\in I$ with $m_{ij} \neq 1$. We call the
substitution $\bi(i,j)_{m_{ij}}\bj \rightarrow
\bi(j,i)_{m_{ji}}\bj$ a {\it braid move}, qualifying it {\it
short} or {\it long} according as $m_{ij}$ equals $2$ or $3$.

Let $w\in W$.  A well-known result of Matsumoto \cite{{\bf 7}} and
Tits \cite{{\bf 9}} states that any reduced expression for $w$ can
be transformed into any other by applying a (possibly empty)
sequence of braid moves.

We say that two words are {\it commutation equivalent} if one can
be transformed into the other by a sequence of short braid moves.
The set of words that are commutation equivalent to a given word
is called the {\it commutation class} of that word.  If the set of
reduced expressions for an element $w\in W$ forms a single
commutation class, then we call $w$ {\it fully commutative},
following \cite{{\bf 8}, \S1}.

Applying a braid relation to a reduced expression corresponds to
applying a permutation to the root sequence of that reduced
expression.  The following proposition makes this more precise.

\proclaim{Proposition 1.2.1 \cite{{\bf 5}, Proposition 3.1.1}} Let
$w \in W$, let $\bi, \bj \in I^*$ and let $i,j,k\in I$.  Denote
the length of $\bj$ by $n$.
\item{\rm (a)}{Assume that $\bi i j \bj$ is a reduced
expression for $w$, and let $\brr = (r_l)$ be the associated root
sequence.
\item{\rm (i)}{Suppose $m_{ij} = 2$, so that
$\bi j i \bj$ is also a reduced expression for $w$. Then the root
sequence $\brr'$ of $\bi j i \bj$ can be obtained from $\brr$ by
interchanging $r_{n+1}$ and $r_{n+2}$, which are mutually
orthogonal relative to $B$.}
\item{\rm (ii)}{If $r_{n+1}$ and $r_{n+2}$ are orthogonal,
then $m_{ij} = 2$.} }
\item{\rm (b)}{Assume that $\bi i j k \bj$ is a reduced
expression for $w$, and let $\brr = (r_l)$ be the associated root
sequence.
\item{\rm (i)}{Suppose $k = i$, so that $m_{ij}=3$ and
$\bi j i j \bj$ is also a reduced expression for $w$. Then the
root sequence $\brr'$ of $\bi j i j \bj$ can be obtained from
$\brr$ by interchanging $r_{n+1}$ and $r_{n+3}$. Furthermore, we
have $r_{n+1} + r_{n+3}=r_{n+2}$.}
\item{\rm (ii)}{If $r_{n+1} + r_{n+3} = r_{n+2}$, then
$k = i \ne j$ and $m_{ij} = 3$.} } \qed\endproclaim

Let $\brr$ and $\brr'$ be as in part (a)(i) (respectively, part
(b)(i)) of Proposition 1.2.1. Employing again the terminology used
above for words, we say that the passage from $\brr$ to $\brr'$ is
obtained by a {\it short braid move} (respectively, {\it long
braid move}). Two root sequences are said to be {\it commutation
equivalent} if one can be transformed into the other by applying a
sequence of short braid moves. The set of root sequences that are
commutation equivalent to a given root sequence is called the {\it
commutation class} of that root sequence.

Let $w\in W$.  The recipe for associating a root sequence to a
reduced expression defines a bijection from the set of reduced
expressions for $w$ to the set of root sequences for $w$, and by
Proposition 1.2.1, this bijection is compatible with the
application of both long and short braid moves.  Hence, by the
result of Matsumoto and Tits cited above, any root sequence for
$w$ can be transformed into any other by applying a sequence of
long and short braid moves. It also follows that there is a
natural bijection between the set of commutation classes of
reduced expressions for $w$ and the set of commutation classes of
root sequences for $w$.

A {\it subword} of a word $i_1i_2\cdots i_n$ (each $i_l\in I$) is
any word of the form $i_pi_{p+1}\cdots i_q$, where $1\leq p \leq q
\leq n$.

\proclaim{Proposition 1.2.2} Let $w \in W$. The following are
equivalent:
\item{\rm (i)}{$w$ is fully commutative.}
\item{\rm (ii)}{$w$ has no inversion triples.}
\item{\rm (iii)}{$w$ has no contractible inversion triples.}
\item{\rm (iv)}{No reduced expression for $w$ contains a subword
of the form $iji$, where $i, j\in I$.}
\endproclaim

\demo{Proof}
\item{(i) $\Rightarrow$ (ii)}
{This is the implication (a) $\Rightarrow$ (c) of \cite{{\bf 3},
Theorem 2.4}.}
\item{(ii) $\Rightarrow$ (iii)}
{This is immediate from the definitions.}
\item{(iii) $\Rightarrow$ (iv)}
If $w$ has a reduced expression with a subword of the form $iji$,
then Proposition 1.2.1 (b)(i) shows that $w$ has a contractible
inversion triple.
\item{(iv) $\Rightarrow$ (i)}
If $w$ is not fully commutative, then there exists a pair of
commutation inequivalent reduced expressions for $w$. It follows
by the result of Matsumoto and Tits mentioned above that $w$ has a
reduced expression to which a braid move can be applied. Thus, $w$
does not satisfy (iv). \qed\enddemo

\head 2. Freely braided elements and commutation classes \endhead

\subhead 2.1 The map $F_w$ \endsubhead

Let $w \in W$. Fix an arbitrary antisymmetric relation $\preceq$
on $\Phi(w)$ with the property that any two roots in $\Phi(w)$ are
comparable relative to $\preceq$. Let ${\Cal C}(w)$ and ${\Cal
I}(w)$ denote the set of commutation classes of root sequences for
$w$ and the set of contractible inversion triples of $w$,
respectively.  We define a map
$$ F_w : {\Cal C}(w) \longrightarrow
\{0,1\}^{{\Cal I}(w)},$$ depending on $\preceq$, as follows. If
${\Cal I}(w)$ is empty, then $\{0,1\}^{{\Cal I}(w)}$ contains just
the empty map, and the set ${\Cal C}(w)$ is also a singleton by
Proposition 1.2.2. Thus, in this situation, $F_w$ is uniquely
determined. Suppose that ${\Cal I}(w)$ is nonempty. Let $C\in
{\Cal C}(w)$ and let $\leq_C$ be the partial ordering of $\Phi(w)$
obtained by taking the transitive closure of the following
relations: $\a < \be$ whenever $\a$ lies to the left of $\be$ in
some root sequence from $C$ and $B(\a, \be)\ne 0$.  Note that
$\leq_C$ is well-defined by \cite{{\bf 5}, Proposition 3.1.5}.
Given any $\{\a,\be, \a+\be \}\in {\Cal I}(w)$, we define
$F_w(C)(\{\a, \be, \a+ \be \})$ to be $0$ if $\a$ and $\be$ are in
the same relative order with respect to $\leq_C$ and $\preceq$,
and otherwise we define $F_w(C)(\{\a, \be, \a+ \be \})$ to be $1$.

The map $F_w$ is injective \cite{{\bf 5}, Theorem 4.1.1}.

It will be convenient to have the following terminology when
determining the surjectivity, or otherwise, of $F_w$. Let $w\in W$
and let ${\Cal T}$ be a subset of ${\Cal I}(w)$. We say that $F_w$
{\it separates} ${\Cal T}$ if every map from ${\Cal T}$ to
$\{0,1\}$ is the restriction of some element of $F_w({\Cal
C}(w))$. Clearly, if $F_w$ fails to separate some nonempty subset
of ${\Cal I}(w)$, then $F_w$ is not surjective.

\subhead 2.2 First main result \endsubhead

It was shown in \cite{{\bf 5}, Corollary 4.1.2, Corollary 4.2.4}
that every $w\in W$ has at most $2^{N(w)}$ commutation classes,
with equality if $w$ is freely braided. The following theorem
shows that equality is achieved only if $w$ is freely braided. For
the special case where $W$ is of type $A$, this was already done
in \cite{{\bf 5}, Theorem 5.2.1} using an ad hoc argument.

\proclaim{Theorem 2.2.1} If $w \in W$ has $2^{N(w)}$ commutation
classes, then $w$ is freely braided.
\endproclaim

\demo{Proof} Let $w$ be a non-freely-braided element. Since $F_w$
is injective, it suffices to prove that $F_w$ is not surjective.
Let $\a$ be a root belonging to at least two contractible
inversion triples of $w$, and assume that the height of $\a$ is
maximal with respect to this property. Let $T$ and $T'$ be
distinct contractible inversion triples of $w$ containing $\a$.
Note that $|T\cap T'|=1$ (this follows easily from \cite{{\bf 5},
Remark 2.2.2} and the contractibility of the triples). By
symmetry, there are three cases to consider.

\medskip

\noindent {\bf Case 1:} $T = \{\a,\be, \a+\be \}$ and $T' = \{\a,
\g, \a+\g \}$.

By \cite{{\bf 5}, Remark 2.2.2} and the contractibility of $T$,
$w$ has a root sequence of the form
$$ (\dots,\a,\a+\be,\be,\dots,\a+\g,\dots,\g,\dots)
$$ or $$
(\dots,\g,\dots,\a+\g,\dots,\a,\a+\be,\be,\dots).
$$ We assume the existence of a sequence of the former type, the
argument for the latter being similar.  By our choice of $\a$, the
roots $\a+\be$ and $\a+\g$ cannot belong to the same contractible
inversion triple of $w$. Suppose that $\a+\be$ is not orthogonal
to $\a+\g$. Then, by \cite{{\bf 5}, Proposition 3.2.2}, $\a +\be$
lies to the left of $\a+ \g$ in every root sequence for $w$, and
so it is impossible for $\a$ to be at the same time to the left of
$\a+\be$ and to the right of $\a+\g$. Thus, $F_w$ does not
separate $\{T,T'\}$.

Suppose instead that $\a+\be$ is orthogonal to $\a+\g$. Then
$\a+\be$ and $\g$ are not mutually orthogonal, since $\a+\be$ is
not orthogonal to $\a$.  Furthermore, $\a+\be$ and $\g$ cannot
belong to the same contractible inversion triple of $w$, by our
choice of $\a$. It follows (again by \cite{{\bf 5}, Proposition
3.2.2}) that $\a+\be$ lies to the left of $\g$ in every root
sequence for $w$. This means that in any root sequence for $w$ in
which $\a$ lies to the left of $\a+\be$, the root $\g$ necessarily
lies to the right of $\a+\g$ (otherwise, $\g$ lies to the left of
$\a+\g$, which must then be to the left of $\a$, which in turn is
to the left of $\a+\be$, a contradiction). Again, $F_w$ does not
separate $\{T,T'\}$.

\medskip

\noindent {\bf Case 2:} $T = \{\a, \be, \a-\be\}$ and $T' = \{\a,
\g,  \a-\g\}$.

Here, we may assume without loss of generality that $w$ has a root
sequence of the form $$
(\dots,\g,\dots,\be,\a,\a-\be,\dots,\a-\g,\dots).
$$

Note that $\g$ cannot be orthogonal to both $\be$ and $\a-\be$, or
it would be orthogonal to their sum. We may assume that $\g$ is
not orthogonal to $\be$. If no contractible inversion triple of
$w$ contains both $\g$ and $\be$, then $\g$ lies to the left of
$\be$ in every root sequence for $w$, and consequently there is no
root sequence for $w$ in which $\be$ lies to the left of $\a$ and
$\g$ lies to the right of $\a$. It follows that $F_w$ does not
separate $\{T,T'\}$.

Suppose instead that there is a contractible inversion triple
$T''$ of $w$ that contains $\g$ and $\be$. We claim that $F_w$
does not separate $\{T,T',T''\}$.  To see this, observe that if
$C$ is a commutation class relative to which $\g$ lies to the left
of $\be$ (this determines $F_w(C)(T'')$), and $\be$ lies to the
left of $\a$ (this determines $F_w(C)(T))$, then $\g$ lies to the
left of $\a$ (so that $F_w(C)(T')$ is also determined). (These
conditions are well-defined by \cite{{\bf 5}, Proposition 3.1.5}.)

\medskip

\noindent {\bf Case 3:} $T = \{\a, \be, \a-\be\}$ and $T' = \{\a,
\g, \a+\g \}$.

In this situation, $w$ has a root sequence of the form
$$ (\dots,\be,\a,\a-\be,\dots,\a+\g,\dots,\g,\dots)
$$ or $$
(\dots,\g,\dots,\a+\g,\dots,\be,\a,\a-\be,\dots).
$$  We deal with the former sequence, the analysis of the latter
being similar.

Note that $\a + \g$ cannot be orthogonal to both $\be$ and
$\a-\be$, or it would be orthogonal to their sum.  We may assume
that $\a+\g$ is not orthogonal to $\be$.  By our choice of $\a$,
the roots $\a+\g$ and $\be$ cannot belong to the same contractible
inversion triple of $w$. Hence, $\be$ lies to the left of $\a+\g$
in every root sequence for $w$. It follows that there is no root
sequence for $w$ in which $\a$ lies to the right of $\a + \g$ and
$\be$ lies to the right of $\a$.  This means that $F_w$ does not
separate $\{T,T'\}$.

We conclude that $F_w$ is not surjective. \qed\enddemo

\proclaim{Corollary 2.2.2} Every $w\in W$ has at most $2^{N(w)}$
commutation classes, with equality if and only if $w$ is freely
braided. \qed\endproclaim

\head 3. Free braidedness and full commutativity \endhead

The goal of this section, achieved by Theorem 3.3.3, is to prove
that $W$ has finitely many freely braided elements if and only if
it has finitely many fully commutative elements.

\subhead 3.1 Reduced expressions for freely braided elements
\endsubhead

For the purposes of the proof of Theorem 3.3.3, we wish to have a
clearer picture on the nature of reduced expressions for freely
braided elements.

\definition{Definition 3.1.1}
Let $\bi$ be a word in $I^*$ and suppose that $\bi$ can be written
as $\bu_0\bb_1 \bu_1\bb_2\bu_2 \cdots \bb_p\bu_p$, where each
$\bb_l$ is of the form $iji$ for some $i,j\in I$ with $m_{ij}=3$.
Then we call $\bb_1,\bb_2, \dots ,\bb_p$ a {\it braid sequence}
for $\bi$. If, furthermore, $\bi$ is reduced and $w=\phi(\bi)$ is
freely braided, then we say that $\bi$ is {\it contracted}
provided there exists a braid sequence for $\bi$ with $p=N(w)$
terms.
\enddefinition

\definition{Definition 3.1.2}
Let $\bi \in I^*$.  A word $\bj \in I^*$ is said to be {\it close}
to $\bi$ if there is a (possibly empty) braid sequence $\bb_1,
\bb_2, \dots, \bb_p$ for $\bi$ such that $\bj$ is the word
obtained by applying a long braid move to each of the $\bb_l$. We
also say that $\bj$ is close to $\bi$ {\it via the sequence}
$\bb_1, \bb_2, \dots, \bb_p$.
\enddefinition

Note that if two words are close to one another, then they
represent the same element of $W$.  Note also that any expression
that is close to a contracted reduced expression is itself
contracted.

We say that $i,j \in I$ are $m${\it -commuting}, or simply {\it
commuting}, if $m_{ij} \neq 3$.

\proclaim{Proposition 3.1.3} Let $w\in W$ be freely braided.
\item{\rm (i)}{There exists a contracted reduced expression
$\bi$ for $w$.}
\item{\rm (ii)}{The reduced expressions close to $\bi$, which are
also contracted reduced expressions, form an irredundantly
described set of commutation class representatives for $w$.}
\item{\rm (iii)}{Any contracted reduced expression for $w$ has
a unique braid sequence with $N(w)$ terms.}
\endproclaim

\demo{Proof} By \cite{{\bf 5}, Theorem 4.2.3}, there is a root
sequence $\brr$ for $w$ such that the roots in any given
contractible inversion triple of $w$ appear consecutively in
$\brr$.  Part (i) follows by applying Proposition 1.2.1 (b)(ii):
take $\bi$ to be the reduced expression corresponding to $\brr$,
and note that the $N(w)$ contractible inversion triples correspond
to a braid sequence for $\bi$ with $N(w)$ terms.  In view of
Proposition 1.2.1 (b)(i), the expression $\bi$, or any contracted
reduced expression for $w$, has at most one braid sequence with
$N(w)$ terms.  This proves (iii).

If $\bi',\bi'' \in I^*$ are distinct and close to $\bi$, then
$\bi'$ is not commutation equivalent to $\bi''$ (to see this,
observe that the sequence of occurrences of any pair of
non-$m$-commuting letters in a word is an invariant of the
commutation class of that word). Therefore, since there are
$2^{N(w)}$ expressions close to $\bi$, and since $w$ has exactly
$2^{N(w)}$ commutation classes (by Corollary 2.2.2), part (ii) is
proved. \qed\enddemo

\remark{Remark 3.1.4} Let $w\in W$.  From the proof of Proposition
3.1.3 (ii), we see that if a reduced expression for $w$ has a
braid sequence with $p$ terms, then $w$ has at least $2^p$
commutation classes.
\endremark

\definition{Definition 3.1.5}
Let $\bi$ be a contracted reduced expression for a freely braided
element $w \in W$, and write $$ \bi = \bu_1 \bb_1 \bu_2 \bb_2
\cdots \bu_{N(w)} \bb_{N(w)} \bq ,$$ where
$\bb_1,\bb_2,\dots,\bb_{N(w)}$ is the unique braid sequence for
$\bi$ with $N(w)$ terms.  If $N(w)>0$, then we define
$D(\bi)=D^1(\bi)$ to be the word in $I^*$ obtained from $\bi$ by
deleting the rightmost letter in $\bb_{N(w)}$.  We do not define
$D(\bi)$ if $N(w)=0$. By induction, we write $D^n(\bi)$ for
$D(D^{n-1}(\bi))$ if $n > 1$ and the composition is defined.  We
agree that $D^0(\bi) = \bi$, regardless of the value of $N(w)$.
\enddefinition

Our strategy for the proof of Theorem 3.3.3 will be to argue that
if $\bi$ is a contracted reduced expression for a freely braided
element $w \in W$, then $D^{N(w)}(\bi)$ is a well-defined reduced
expression for some fully commutative element. This will require
several intermediate steps. One of these is the following
technical lemma.

\proclaim{Lemma 3.1.6} Maintain the notation of Definition 3.1.5.
Suppose that $D(\bi)$ is a reduced expression for a freely braided
element $w' \in W$ with $N(w') = N(w) - 1$. Then any expression
close to $D(\bi)$ is of the form $D(\bj)$, where $\bj$ is a
reduced expression for $w$ that is close to $\bi$ via a braid
sequence not involving $\bb_{N(w)}$.
\endproclaim

\demo{Proof} The hypotheses on $D(\bi)$ imply that it is a
contracted reduced expression for $w'$, and that $\bb_1, \bb_2,
\dots, \bb_{N(w) - 1}$ is a braid sequence for $D(\bi)$ with
$N(w')$ terms.  The conclusion follows. \qed\enddemo

\subhead 3.2 Freely braided elements and the $N$-statistic \endsubhead

The following lemma describes what happens when one goes up in the
weak Bruhat order from a freely braided element.

\proclaim{Lemma 3.2.1} Suppose that $w\in W$ is freely braided,
and that $\ell(ws_i) > \ell(w)$ for some $i \in I$. As usual, we
denote the simple root corresponding to $i$ by $\a_i$.
\item{\rm (a)}{Assume that $\a_i$ does not lie in any
contractible inversion triple of $ws_i$.
\item{\rm (i)}{The roots occurring before $\a_i$ in any
root sequence for $ws_i$ are orthogonal to $\a_i$.}
\item{\rm (ii)}{The contractible inversion triples of $ws_i$ are
precisely those of the form $s_i(T)$, where $T$ is a contractible
inversion triple of $w$.}
\item{\rm (iii)}{The element $ws_i$ is freely braided and
$N(ws_i) = N(w)$.} }
\item{\rm (b)}{Assume that $\a_i$ lies in some contractible
inversion triple of $ws_i$.
\item{\rm (i)}{There is a reduced expression $\bi$ for $w$ of the
form $\bu i j \bv, $ where $j \in I$ does not commute with $i$ and
where each letter in $\bv$ commutes with $i$.}
\item{\rm (ii)}{Any reduced expression for $w$ that is commutation
equivalent to the reduced expression $\bi$ in {\rm (i)} must be of
the form $ \bu' i \bv_1 j \bv_2 ,$ where each letter in $\bv_1$
and each letter in $\bv_2$ commutes with $i$.} }
\endproclaim

\demo{Note} We do not require above that $ws_i$ be freely braided.
\enddemo

\demo{Proof}
We first prove (a).

Since there is a reduced expression for $ws_i$ in which $s_i$
appears last, there is a root sequence $\brr$ for $ws_i$ in which
$\a_i$ is the first root. Let $\brr'$ be an arbitrary root
sequence for $ws_i$.  By the discussion following Proposition
1.2.1, $\brr'$ may be obtained from $\brr$ by applying a sequence
of braid moves. Since none of these braid moves can be a long
braid move involving $\a_i$, we find that all the roots occurring
before $\a_i$ in $\brr'$ are orthogonal to $\a_i$. This proves
(i).

Suppose that $T$ is a contractible inversion triple of $w$, and
let $\brr_0$ be a root sequence for $w$ in which the elements of
$T$ appear consecutively.  Since $(\a_i, s_i(\brr_0))$ is a root
sequence for $ws_i$ (recall the definition of root sequence in
\S1.1), it follows that $s_i(T)$ is a contractible inversion
triple of $ws_i$.

Conversely, suppose that $T'$ is a contractible inversion triple
of $ws_i$, and let $\brr_1$ be a root sequence for $ws_i$ in which
the elements of $T'$ appear consecutively.  By hypothesis, $\a_i$
does not appear in $T'$, and by (i), the elements appearing before
$\a_i$ in $\brr_1$ are orthogonal to $\a_i$. Hence, we may apply
short braid moves if necessary to obtain a root sequence $\brr'_1$
for $ws_i$ in which $\a_i$ appears first and in which the elements
of $T'$ still appear consecutively. Now, $\brr'_1$ is of the form
$(\a_i, s_i(\brr''_1))$, where $\brr''_1$ is a root sequence for
$w$ in which the elements of $s_i(T')$ occur consecutively. This
proves (ii).

By (ii), we have $N(ws_i) = N(w)$.  Let $\brr_2$ be a root
sequence for $w$ of the form specified in \cite{{\bf 5}, Theorem
4.2.3}. Using (ii) again, we see that the contractible inversion
triples of $ws_i$, the roots in each of which appear consecutively
in the root sequence $(\a_i, s_i(\brr_2))$, are pairwise disjoint.
Hence, $ws_i$ is freely braided, and (iii) is proved.

We turn to (b).

Let $\brr$ be a root sequence for $ws_i$ in which $\a_i$ appears
first, and consider a sequence of braid moves of minimal length
subject to the condition that applying the sequence to $\brr$
results in a root sequence $\brr'$ in which the elements of some
contractible inversion triple containing $\a_i$ appear
consecutively. Denote by $T$ the contractible inversion triple
that contains $\a_i$ and is consecutive in $\brr'$. By the
minimality assumption, none of the braid moves in the above
sequence is a long braid move involving $\a_i$. Hence, every root
occurring before $\a_i$ in $\brr'$ is orthogonal to $\a_i$, and we
may therefore apply a sequence of short braid moves to $\brr'$ to
obtain a root sequence $\brr''$ in which $\a_i$ appears first and
in which the other elements of $T$ appear consecutively. By
Proposition 1.2.1, the sequence $\brr''$ corresponds to a reduced
expression for $ws_i$ of the form $\bu i j \bv i$, where all of
the letters in $\bv$ commute with $i$.  Deleting the rightmost
$i$, we obtain a reduced expression for $w$ of the required form,
thus proving (i).

To prove (ii), we note that the reduced expression obtained in (i)
is of the stated form, taking $\bu'=\bu$, $\bv_1 = \emptyset$ and
$\bv_2 = \bv$. The result follows, once we observe that applying a
short braid move to a reduced expression of the form given in (ii)
produces another expression of the same form. \qed\enddemo

The next result describes what happens when one goes down in the
weak Bruhat order from a freely braided element.

\proclaim{Lemma 3.2.2} Suppose that $w\in W$ is freely braided and
that $i \in I$ satisfies $\ell(ws_i) < \ell(w)$.   Then $ws_i$ is
freely braided, and we have $$ N(ws_i) = \cases N(w) - 1 & \text{
if } \a_i \text{ lies in a contractible inversion triple of
$w$},\cr N(w) & \text{ otherwise.}\cr
\endcases
$$
\endproclaim

\demo{Proof}  If $\brr$ is a root sequence for $ws_i$, then
$(\a_i, s_i(\brr))$ is a root sequence for $w$. Therefore, every
contractible inversion triple $T$ of $ws_i$ yields a contractible
inversion triple $s_i(T)$ of $w$. This gives $N(w) \geq N(ws_i)$.
The inequality is strict if $\a_i$ lies in a contractible
inversion triple of $w$.

By the Exchange Condition for Coxeter groups (see \cite{{\bf 6},
\S 5.8}), $w$ has a reduced expression $\bi$ ending with $i$, and
by Proposition 3.1.3 (ii), there is a contracted reduced
expression $\bj$ for $w$ that is commutation equivalent to $\bi$.
Write $\bj = \bv_1 i \bv_2$, where each letter in $\bv_2$ commutes
with $i$. Suppose that $\a_i$ lies in a contractible inversion
triple of $w$.  Then $N(ws_i) \leq N(w) - 1$ by the first
paragraph.  On the other hand, it is clear that $\bv_1 \bv_2$, a
reduced expression for $ws_i$, has a braid sequence with $N(w)-1$
terms; hence, $N(ws_i) \geq N(w) - 1$ by Proposition 1.2.1 (b)(i).
Further, by Remark 3.1.4, $ws_i$ has at least $2^{N(ws_i)}$
commutation classes.  It now follows from Corollary 2.2.2 that
$ws_i$ is freely braided.

Suppose instead that $\a_i$ does not belong to a contractible
inversion triple of $w$. Then, by Proposition 1.2.1 (b)(i), $\bv_1
\bv_2$ has a braid sequence with $N(w)$ terms.  It follows by the
same proposition together with the first paragraph that $N(ws_i) =
N(w)$. As above, we find that $ws_i$ has at least $2^{N(ws_i)}$
commutation classes, and so is freely braided. \qed\enddemo

\remark{Remark 3.2.3} If $w\in W$ is freely braided and
$\ell(ws_i) > \ell(w)$ with $ws_i$ non-freely-braided, it may
happen that $N(ws_i) > N(w) + 1$. For example, if $W$ is of type
$A_3$ and $w = s_2 s_1 s_3 s_2$ (using the obvious indexing), then
$N(w) = 0$ but $N(w s_3) = 2$.
\endremark

\subhead 3.3 Groups with finitely many freely braided elements \endsubhead

The following lemma is a crucial ingredient in the proof of
Theorem 3.3.3.

\proclaim{Lemma 3.3.1} Let $\bi$ be a contracted reduced
expression for a freely braided element $w \in W$ with $N(w)
>0$. Then $D(\bi)$ is a contracted reduced expression for a freely
braided element $w'$ with $N(w') = N(w) - 1$.
\endproclaim

\demo{Proof} We start by writing $$ \bi = \bc_1 \bb_1 \bc_2 \bb_2
\cdots \bc_{N(w)} \bb_{N(w)} \bq,$$ where $\bb_1,\bb_2,\dots
,\bb_{N(w)}$ is the unique braid sequence for $\bi$ with $N(w)$
terms.  Let $\bi_q$ be the expression obtained from $\bi$ by
deleting from $\bq$ all but its first $q$ letters. We thus have
$\bi=\bi_l$, where $l$ is the length of $\bq$. If $\brr$ is the
root sequence of $\bi$, then the first $l$ roots in $\brr$ are not
involved in any contractible inversion triple of $w$, and it
follows from repeated applications of Lemma 3.2.2 that for each
$q$, the reduced expression $\bi_q$ represents a freely braided
element $y_q$ with $N(y_q) = N(w)$.  Moreover, if $q > 0$, then
the first root in the root sequence of $\bi_q$ does not lie in any
contractible inversion triple of $y_q$.

We will prove by induction on $q$ that $D(\bi_q)$ is a reduced
expression for a freely braided element $w_q$ with $N(w_q) = N(w)
- 1$. (These properties imply that $D(\bi_q)$ is contracted.)
Denote the letter that is deleted from $\bi$ to form $D(\bi)$ by
$j$.

\medskip

{\noindent \bf Base case: $q = 0$.}

In this case, $D(\bi_0)$ is obtained from $\bi_0$ by removing the
last letter, $j$.  It is clear that $D(\bi_0)$ is a reduced
expression for some group element $w_0$, which is freely braided
by Lemma 3.2.2.  The first three roots in the root sequence of
$\bi_0$ comprise a contractible inversion triple containing
$\a_j$. Hence, by Lemma 3.2.2 again, $N(w_0)=N(y_0)-1$, and the
latter equals $N(w)-1$ by the first paragraph.

\medskip

{\noindent \bf Inductive step: proof that $D(\bi_q)$ is reduced.}

Suppose that the statement is true for all $q$ with $0\leq q \leq
k < l$. Let $q=k+1$, and let $i$ be the $(k+1)$-st factor of
$\bq$. By the inductive hypothesis, $D(\bi_k)$ is a contracted
reduced expression for a freely braided element $w_k$ with $N(w_k)
= N(w) - 1$. Assume toward a contradiction that $D(\bi_{k+1})$ is
not reduced. Then by the Exchange Condition for Coxeter groups,
there is a reduced expression $\bi'$ for $w_k$ ending in $i$.
According to Proposition 3.1.3 (ii), there is a unique contracted
reduced expression $\bi''$ for $w_k$ that is both close to
$D(\bi_k)$ and commutation equivalent to $\bi'$. We may write
$$ \bi'' = \bp'' i \bc'' ,$$ where all of the letters in
$\bc''$ commute with $i$.

By Lemma 3.1.6, $\bi''$ must be of the form $D(\bj)$, where $\bj$
is close to $\bi_k$. The expression $\bj$ is thus obtainable from
$\bi''$ by inserting the letter $j$ at some point into the word.
Since $\bi_{k+1}$ is reduced, this insertion must take place to
the right of the indicated occurrence of $i$ in $\bi''$.
Therefore, $\bj i$, which is a reduced expression for $y_{k+1}$,
is of the form $$\bp'' i \bc_1'' j \bc_2'' i ,$$ where each letter
in $\bc_1''$ and each letter in $\bc_2''$ commutes with $i$.
Applying short braid moves if necessary, we obtain
$$\bp'' \bc_1'' i j i \bc_2'',$$ and it follows from Proposition
1.2.1 (b) that the first root in the root sequence of $\bj i$
belongs to a contractible inversion triple. Now, $\bj$
(respectively, $\bj i$) is a reduced expression for $y_k$
(respectively, $y_{k+1}$), and $N(y_k) = N(y_{k+1}) = N(w)$. This
contradicts Lemma 3.2.2, taking $w = y_{k+1}$.

We conclude that $D(\bi_{k+1})$ is reduced.

\medskip

{\noindent \bf Inductive step: proof that $w_q$ is freely braided
and $N(w_q) = N(w) - 1$.}

By the above, we have $w_{k+1}=w_ks_i$ with $\ell(w_ks_i) >
\ell(w_k)$. If $\a_i$ does not lie in any contractible inversion
triple of $w_ks_i$, then the inductive step follows from Lemma
3.2.1 (a)(iii) (with $w_k$ playing the role of $w$).  Assume
instead that we are in case (b) of Lemma 3.2.1, which is the only
alternative.

By Lemma 3.2.1 (b)(i), the element $w_k$ has a reduced expression
of the form $\bu i i' \bv,$ where $i'\in I$ does not commute with
$i$ and every letter in $\bv$ commutes with $i$. Recall that
$D(\bi_k)$ is contracted by the inductive hypothesis. Hence, by
Proposition 3.1.3 (ii), there is a contracted reduced expression
$\bd_k$ for $w_k$ that is both close to $D(\bi_k)$ and commutation
equivalent to $\bu i i' \bv$. According to Lemma 3.2.1 (b)(ii), we
have $\bd_k = \bu' i \bv_1 i' \bv_2 ,$ where every letter in
$\bv_1$ and every letter in $\bv_2$ commutes with $i$. Because
$\bd_k$ is close to $D(\bi_k)$ and $N(w_k)=N(y_k)-1$, Lemma 3.1.6
implies that $\bd_k$ is of the form $D(\bi'_k)$, where $\bi'_k$ is
a reduced expression for $y_k$ that is close to $\bi_k$ by a
sequence of braid relations not involving $\bb_{N(w)}$. There is
no loss in generality in assuming that $\bi_k$ is equal to
$\bi'_k$, so we will do this in order to make the arguments
clearer.

Let $\bd_{k+1} = \bd_ki = \bu' i \bv_1 i' \bv_2 i$. Since $\bd_k =
D(\bi_k)$, we have $\bd_{k+1} = D(\bi_{k+1})$. Hence, by
appropriately inserting $j$ in $\bd_{k+1}$, we obtain $\bi_{k+1}$.
The insertion must take place immediately to the right of a
subword of $\bd_{k+1}$ of the form $jj'$, where $j'\in I$ does not
commute with $j$.

We know from the first paragraph of the proof that $\a_i$ does not
lie in a contractible inversion triple of $y_{k+1}$.  The only way
this can happen is if the letter $j$ is inserted in $\bd_{k+1}$
somewhere between the two indicated occurrences of $i$, and if $j$
does not commute with $i$.  Since the letter sitting two places to
the left of the insertion site is also an occurrence of $j$, the
latter occurrence of $j$ must either occur to the left of the
leftmost indicated occurrence of $i$ in $\bd_{k+1}$, or must be
the indicated occurrence of $i'$ in $\bd_{k+1}$.  We consider each
of these two cases in turn.

In the first case, $\bi_{k+1}$ can be written as
$$ \bu'' j i j \bv_1 i' \bv_2 i ,$$ where $\bv_1$ and
$\bv_2$ are as before and $\bu''j=\bu'$. Applying a long braid
move, we obtain the following reduced expression for $y_{k+1}$:
$$ \bu'' i j i \bv_1 i' \bv_2 i .$$  This contradicts the
fact (mentioned in the first paragraph of the proof) that
$y_{k+1}$ is freely braided, because Proposition 1.2.1 can now be
used to show that the middle of the three indicated occurrences of
$i$ corresponds to a root that lies in two different contractible
inversion triples.

In the second case, $i'=j$ and $\bi_{k+1}$ can be written as
$$ \bu' i \bv_1 j j' j \bv'_2 i ,$$ where $\bv_1$ is as
above and $j'\bv'_2=\bv_2$, meaning that $i$ commutes with $j'$.
Applying a long braid move and commutations, we obtain
$$ \bu' \bv_1 i j' j j' i \bv'_2 .$$ This again leads to a
contradiction because the indicated occurrence of $j$ corresponds
to a root that belongs to two different contractible inversion
triples of $y_{k+1}$.

We have completed the inductive step by showing that case (b) of
Lemma 3.2.1 cannot occur. \qed\enddemo

\proclaim{Corollary 3.3.2} Let $\bi$ be a contracted reduced
expression for a freely braided element $w \in W$, and maintain
the notation of Definition 3.1.5, so that $$ \bi = \bu_1 \bb_1
\bu_2 \bb_2 \cdots \bu_{N(w)} \bb_{N(w)} \bq .$$ Then the
expression obtained by omitting the rightmost letter in each of
the words $\bb_l$ is a reduced expression for a fully commutative
element.
\endproclaim

\demo{Proof} By applying Lemma 3.3.1 $N(w)$ times, we find that
$D^{N(w)}(\bi)$, which is the expression described in the
conclusion, is a reduced expression for a freely braided element
$w'$ with $N(w') = 0$.  By Proposition 1.2.2, $w'$ is fully
commutative. \qed\enddemo

\proclaim{Theorem 3.3.3} A simply laced Coxeter group $W$ has
finitely many freely braided elements if and only if it has
finitely many fully commutative elements.
\endproclaim

\demo{Proof} By Proposition 1.2.2, any fully commutative $w \in W$
satisfies $N(w) = 0$, and so is freely braided for vacuous
reasons.  This proves the ``only if" part of the theorem.

Conversely, suppose that $W$ has finitely many fully commutative
elements. By Corollary 3.3.2, there is a map from the set of
contracted reduced expressions for freely braided elements of $W$
to the set of reduced expressions for fully commutative elements
of $W$, given by $$ \bi \mapsto D^{N(\phi(\bi))}(\bi) .$$  Any
element in the fibre over $D^{N(\phi(\bi))}(\bi)$  can be
recovered from $D^{N(\phi(\bi))}(\bi)$ by making appropriate
insertions of generators after certain subwords of the form $ij$,
where $i,j\in I$ are noncommuting.  Hence, the fibres of the above
map are all finite. Since there are finitely many fully
commutative elements and each of these has finitely many reduced
expressions, there are finitely many fibres. It follows that $W$
has finitely many freely braided elements. \qed\enddemo

Independently of one another, Graham \cite{{\bf 4}} and Stembridge
\cite{{\bf 8}} have classified the Coxeter groups with finitely
many fully commutative elements. The classification has also been
worked out in the simply laced case by Fan \cite{{\bf 2}}. It
turns out that a simply laced Coxeter group has finitely many
fully commutative elements if and only if each connected component
of its Coxeter graph is of type $A_n$, $D_n$ or $E_n$ for
arbitrary $n$. In particular, Coxeter groups of type $E_n$ for $n
> 8$ have finitely many fully commutative elements, although the
groups are infinite. This classification carries over for freely
braided elements, by the above theorem.

\leftheadtext{}
\rightheadtext{}

\end